\UseRawInputEncoding
\documentclass[12pt,reqno]{amsart}
\usepackage{amssymb,amsmath,graphicx,
	amsfonts,euscript}
\usepackage{color}
\usepackage{mathrsfs}

\newcommand{\beq}{\begin{equation}}
\newcommand{\eeq}{\end{equation}}
\newcommand{\ben}{\begin{eqnarray}}
\newcommand{\een}{\end{eqnarray}}
\newcommand{\beno}{\begin{eqnarray*}}
\newcommand{\eeno}{\end{eqnarray*}}

\newtheorem{Theorem}{Theorem}[section]
\newtheorem{Definition}[Theorem]{Definition}
\newtheorem{Proposition}[Theorem]{Proposition}
\newtheorem{Lemma}[Theorem]{Lemma}

\newtheorem{Remark}[Theorem]{Remark}

\begin{document}

\title[Regularity Criteria Without Pressure for the Navier-Stokes Equations]{Remarks on Interior Regularity Criteria Without Pressure for the Navier-Stokes Equations}

\author{Shuai Li}
\address{School of Mathematical Sciences, Dalian University of Technology, Dalian, 116024, PR China}
\email{leeshy@mail.dlut.edu.cn}

\author{Wendong Wang}
\address{School of Mathematical Sciences, Dalian University of Technology, Dalian, 116024, PR China}
\email{wendong@dlut.edu.cn}

\author{Daoguo Zhou}
\address{School of Mathematics, Hangzhou Normal University, Hangzhou, 311121, PR China}
\email{daoguozhou@qq.com}


\date{\today}


\begin{abstract}
In this note we investigate interior regularity criteria for suitable
	weak solutions to the 3D Naiver-Stokes equations, and obtain the solutions are regular in the interior if the $L^p_tL_x^q(Q_1)$ norm of the velocity is sufficiently small, where $1\leq \frac{2}{p}+\frac{3}{q}<2$ and $2\leq p\leq \infty$. It improves the recent result of $p,q>2 $ by Kwon \cite{Kwon} (J. Differential Equations 357 (2023), 1--31.), and also generalizes Chae-Wolf's $L_t^\infty L_x^{\frac32+}$ criterion \cite{CW2017} (Arch. Ration. Mech. Anal. 225 (2017), no. 1, 549--572.).

	\end{abstract}
	
	\maketitle
	
	
	\section{Introduction}
Consider the 3D Navier-Stokes equations describing viscous incompressible fluid in $ \mathbb{R}^3\times(0,T)$:
\begin{eqnarray}\label{NSE}
 \left\{
    \begin{array}{llll}
    \displaystyle \partial_t u- \Delta u+ u\cdot \nabla u+ \nabla \pi = 0,\\
   \text{div} \; u= 0
    \end{array}
 \right.
\end{eqnarray}
with a smooth and rapidly decaying solenoidal initial vector field $u(x,0)=u_0(x)$ in $\mathbb{R}^3$. Here $u(x,t)$ denotes the velocity of the fluid and the scalar function $\pi(x,t)$ denotes the pressure.

 	
	
	In a seminal paper \cite{Leray34}, Leray proved the global existence of weak solutions with finite energy to the Navier-Stokes equations in three dimensions. See also the global existence of weak solution in a bounded domain by Hopf  \cite{Hopf}. However, the regularity of weak solutions is still an outstanding open problem in mathematical fluid mechanics. One type of condition ensuring regularity is that
\ben\label{ine:GRC}
\|u\|_{L^p((0,T);L^q(\mathbb{R}^3))} < +\infty, \quad \frac 2p + \frac 3 q = 1, \quad q \in [3,+\infty],
\een
and we refer to Ladyzenskaja \cite{Ladyzenskaja}, Prodi \cite{Prodi}, Serrin \cite{Serrin62} Struwe \cite{Struwe} and the references therein.
The endpoint case of $p=\infty, q = 3$ is highly nontrivial, which was resolved by Escauriaza-Seregin-\v{S}ver\'{a}k in \cite{ESS}.

In a series of papers \cite{Scheffer76,Scheffer77}, Scheffer began the partial regularity theory of the Navier-Stoeks equations. Caffarelli, Kohn and Nirenberg \cite{CKN}  improved the results of Scheffer by proving that the set $\mathcal{S}$ of possible interior singular points of a suitable weak solution is of one-dimensional parabolic Hausdorff measure zero, i.e. $\mathcal{P}^1(\mathcal{S})=0$, which rests on the following two $\varepsilon-$regularity criteria for suitable weak solutions to (\ref{NSE}). There is an absolute constant $\varepsilon>0$ such that $u$ is regular at $(0,0)$ if one of the following conditions holds:
	\begin{itemize}
\item[]
	\begin{align}\label{ine:one scale}
	\|u\|_{L^3(Q_1)} + \|u \pi\|_{L^1(Q_1)} + \|\pi\|_{L^\frac54 L^1(Q_1)} \leq \varepsilon,
	\end{align}
	\item[]
	\begin{align}\label{ine:infinity scale}
	\limsup\limits_{r\rightarrow 0^+} r^{-1} \int_{Q_{r}} |\nabla u(y,s)|^2 dyds \leq \varepsilon.
	\end{align}
\end{itemize}
Recall  that suitable weak solutions are defined  as follow:
\begin{Definition}
	\label{sws} Let $\Omega$ be a domain in $\mathbb R^3$ and let $Q_T:=\Omega\times (-T,0)$. We say that a pair of functions $(u,\pi)$ is a suitable weak solution to the Navier-Stokes equations in $Q_T$ if the following conditions are fulfilled:
	
	(i) $u\in L_{loc}^\infty(-T,0;L^2_{loc}(\Omega))\cap L^{2}(-T,0;W^{1,2}_{loc}(\Omega)),\quad \pi\in L^{\frac 32}(-T,0;L^{\frac 32}_{loc}(\Omega))$;

	(ii) $u$ and $\pi$ satisfy the Navier-Stokes equations (\ref{NSE})
	in $Q_T$ in the sense of distributions;
	
	(iii) for $Q(z_0,R)\subset Q_T$, the local energy inequality
	\beno&&\int\limits_{B(x_0,R)}\varphi|u(x,t)|^2dx+2\int\limits^t_{t_0-R^2}\int\limits_{B(x_0,R)}\varphi|\nabla u|^2dxd\tau\\
	&\leq& \int\limits^t_{t_0-R^2}\int\limits_{B(x_0,R)}\bigg (|u|^2(\partial_t\varphi+\Delta \varphi)+u\cdot\nabla \varphi(|u|^2+2\pi)\bigg )dxd\tau\eeno
holds for a.e. $t\in (t_0-R^2,t_0)$ and for all non-negative test functions $\varphi\in C^\infty_0(B(x_0,R) \times (t_0-R^2,t_0+R^2))$.	
\end{Definition}
There are many works devoted to improving local regularity criteria \eqref{ine:one scale} and \eqref{ine:infinity scale}. 
For example, we refer to Tian-Xin \cite{TX}, Ladyzhenskaya-Seregin \cite{LS}, Seregin \cite{Se}, Gustafson-Kang-Tsai \cite{GKT} and the references therein.  In \cite{GKT}, it was proved that $u$ is regular at $(0,0)$, if
	\ben\label{eq: regularity on u}
	\limsup\limits_{r\rightarrow0} r^{1-\frac2p-\frac3q} \|u-\bar{u}\|_{L^p_t L^q_x(Q_r)} \leq \varepsilon, \quad 1\leq \frac2p+\frac3q \leq 2, ~~ 1\leq p,q \leq \infty,
	\een
	or
	\ben\label{eq: regularity on nabla u}
	\limsup\limits_{r\rightarrow0} r^{2-\frac2p-\frac3q} \|\nabla u\|_{L^p_t L^q_x(Q_r)} \leq \varepsilon, \quad 1\leq \frac2p+\frac3q \leq 2, ~~ 1\leq p,q \leq \infty.
	\een
In one scale, Vasseur \cite{ Vasseur} proved $u\in L^\infty(Q_{\frac12})$, if
\beno
\sup_{-1<t<0}\int_{B_1}|u|^2dx+\int_{Q_1}|\nabla u|^2dxdt+\int_{-1}^0\left(\int_{B_1}|\pi|dx\right)^pdt\leq \varepsilon,\quad p>1.
\eeno
Wang-Zhang \cite{WZ14} improved it to
\beno
\sup_{-1<t<0}\int_{B_1}|u|^2dx+\int_{-1}^0\left(\int_{B_1}|u|^4dx\right)^{\frac12}dt+\int_{Q_1}|\pi|dxdt\leq \varepsilon.
\eeno
He, Wang and Zhou  \cite{HWZ} extended \eqref{ine:one scale} as: $u\in L^\infty(Q_{\frac12})$, if
	\beno
	\|u\|_{L^p_t L^q_x(Q_1)} + \|\pi\|_{L^1(Q_1)} \leq \varepsilon, \quad 1\leq \frac2p+\frac3q < 2, \quad 1\leq p,q \leq \infty.
	\eeno
Recently, some local regularity criteria without pressure are established by using the pressure decomposition method which based on Stokes projection.
Wolf in \cite{[Wolf10], [W15]} gave another version of suitable weak solution, which is \emph{}called local suitable weak solution, and proved that $u$ is regular in $Q_\frac {r_0}2$ with some $r_0>0$ if there exists an absolute constant $\varepsilon_0>0$ such that
\ben\label{eq:wolf}
r_0^{-2} \int_{Q_{r_0}} |u|^3 dx dt \leq \varepsilon_0.
\een
Wolf's result were further generalized by Wang, Wu and Zhou in \cite{WWZ} by proving that
\ben\label{eq:wwz}
\|u\|_{L^{q}(Q_1)} \leq \varepsilon_0, \quad q > \frac 52.
\een
Chae-Wolf in \cite{CW2017} also proved the norm of $u\in L_t^\infty L_x^{\frac32+}$ implies the regularity.
Recently, Kwon \cite{Kwon} proved the interior regularity
under the assumption
\beno
\|u\|_{L^r_t L_x^m(Q_2)} \leq \varepsilon_0, \quad  {\rm with} \quad \frac2r + \frac3m < 2, \quad r,m \in (2,+\infty].
\eeno
 by using compactness method for the dissipative weak solution, which was introduced by Duchon-Robert in \cite{DR} firstly. For more results, we refer to \cite{ABP, CML, GP, WWZ} etc.

{\it It's interesting that whether the range of $p,q$ can be relaxed to the range as in \cite{GKT}?}
That is to say,
\beno
 1\leq \frac2p+\frac3q \leq 2, ~~ 1\leq p,q \leq \infty.\eeno
We investigate this issues int this note, and show that $p\geq 2$ is sufficient when $\frac2p+\frac3q = 2-$. Our main result is as follows.
	\begin{Theorem}\label{onescale}Suppose that $(u,\pi)$ is a suitable weak solution to (\ref{NSE}) in $Q_1$.
	For any $(p,q)$ satisfying $1\leq 2/p+3/q<2$, $2\leq p\leq\infty$,
	there exists a positive constant $\varepsilon_0$ such that if
	\[
	\|u\|_{L^{p}L^{q}(Q_{1})}\leq\varepsilon_0,
	\]
	then $u$ is regular at $(0,0)$.
	\end{Theorem}

\begin{Remark} We say that $u$ is regular at a certain point, which means that there exists a neighborhood at this point where u is bounded.
The above theorem improved the result of Kwon \cite{Kwon} by considering the borderline case  $p = 2$ and $q$ is arbitrary, which also generalized Chae-Wolf's $u\in L_t^\infty L_x^{\frac32+}$ in \cite{CW2017} .
\end{Remark}

\begin{Remark}
The restriction $p \geq 2$ in Theorem \ref{onescale} seems to be sharp, which comes from the harmonic part of pressure projection. For example, the term $\int_{Q_\rho} \nabla \pi_h \cdot \nabla \nabla \pi_h \cdot v_B$ (see the term of $M_5$ in (\ref{eq:M5})) in the proof implies $p\geq 2,$
since the harmonic part $\nabla \pi_{h,B}$ and $\nabla \nabla \pi_{h,B}$  can only be  control by $\|u\|_{L_t^pL_x^q}$ due to the iteration.
\end{Remark}

The paper is organized as follows, in Section 2, we introduce some definitions and technical lemmas, especially including the Stokes decomposition of the pressure. In Section 3, we proved a Cacciopolli inequality, which plays an important role in our proof. Theorem {\ref{onescale}} is  proved in Section 4.

Throughout this article, $C$ denotes an absolute constant independent of $u$ and may be different from line to line.



%
%

\section{Preliminaries: some technical lemmas}
\subsection{Local pressure projection}
Let us introduce Wolf's pressure decomposition as in \cite{[Wolf10],[W15]}. For a bounded $C^2$-domain $G\subset \mathbb{R}^n$ and $1<s<\infty$,
 for any $ F\in W^{-1,s}(G)$,  there exists a unique pair $(v,\pi)\in W^{1,s}_{0}\times L^s_0(G)$ which solves the steady Navier-Stokes equations in the weak sense due to  $L^p-$ theory of the steady Stokes system (see, for example, \cite{GSS})
\begin{equation}\label{eq:stokes}
\left\{\begin{array}{llll}
-\Delta v+\nabla \pi=F,\quad &{\rm in}& \quad G,\\
{\rm div }~ v=0, \quad &{\rm in}& \quad G,\\
v=0, \quad &{\rm on}& \quad \partial G,
\end{array}\right.
\end{equation}
where $\pi\in L^s_0(G)$ denotes
\beno
\pi\in L^s(G) \quad {\rm with} \quad \int_G \pi dx=0.
\eeno
Define the operator $E_G$ as follows:
\beno
E_G:W^{-1,s}(G)\rightarrow  W^{-1,s}(G), \quad E_G(F)=\nabla \pi,
\eeno
where $\nabla \pi$ denotes the gradient functional in $W^{-1,s}(G)$ defined by
\beno
<\nabla p,\phi>=-\int_G p \nabla\cdot \phi dx,\quad \phi\in W^{1,s'}_0(G).
\eeno
The operator $E_G$ is bounded from $W^{-1,s}(G)$ into itself with $E_G(\nabla \pi)=\nabla \pi$ for all
$\pi\in L^s_0(G)$, and
\ben\label{eq:bound of wolf}
\|\nabla \pi\|_{L^s(G)}\leq C\|F\|_{W^{-1,s}(G)}.
\een
The norm of $E_G$ depends only on $s$ and the geometric properties of $G$. Specially, the norm of $E_G$ is independent
of $G$, if $G$ is a ball or an annulus, which is due to the scaling properties of the Stokes equation.
Moreover, if $F \in L^s(G)$, by the embedding $L^s(G) \hookrightarrow W^{-1,s}(G)$ and the regularity of elliptic equations, there hold the estimate
\ben\label{eq:bound of wolf'}
\|\pi\|_{L^s(G)}\leq C\|F\|_{L^s(G)}.
\een

\subsection{Local suitable weak solutions}

In this subsection, we define the local suitable weak solution as in \cite{[Wolf10],[W15]}.
\begin{Definition}\label{def:1}
	Let $\Omega$ be a domain in $\mathbb R^3$ and let $Q_T:=\Omega\times (-T,0)$. We say $u$ is a local suitable weak solution to the Navier-Stoeks equations (\ref{NSE}) in $Q_T$ if  \\
(i). $u\in  L^\infty_{\rm loc}(-T,0;L^2_{\rm loc}(\Omega)) \cap L^2_{\rm loc}(-T,0;W^{1,2}_{\rm loc}(\Omega))$;\\
(ii). $u$ is a distributional solution to (\ref{NSE}), i.e. for every $\varphi \in C_c^\infty(\Omega\times (-T,0))$ with $\nabla \cdot \varphi = 0$,
\beno
\int \int_{\Omega\times (-T,0)} -u \cdot \partial_t \varphi - u \otimes u : \nabla \varphi + \nabla u : \nabla \varphi  = 0;
\eeno
(iii). For every ball $B\subset \Omega$ the following energy inequality holds: for almost all $s\in (-T,0)$ and for all non negative $\phi\in C_c^\infty(B\times (-T,0))$, \ben\label{eq:local energy inq} \nonumber
&&\int|v_B(x,s)|^2\phi(x,s)dx+2\int\int|\nabla v_B(x,\tau)|^2\phi(x,\tau) dxd\tau\\ \nonumber
&\leq& \int\int|v_B(x,\tau)|^2(\partial_t\phi+\Delta \phi)dxd\tau+\int\int|v_B|^2(v_B-\nabla \pi_{h,B})\cdot\nabla \phi dxd\tau\\ \nonumber
&&+2\int\int v_B\cdot\nabla \nabla \pi_{h,B} \cdot v_B \phi dxd\tau-2\int\int \nabla \pi_{h,B} \cdot\nabla \nabla \pi_{h,B} \cdot v_B \phi dxd\tau\\
&&+2\int\int (\pi_{1,B}+\pi_{2,B})v_B\cdot \nabla \phi dxd\tau,
\een
with $v_B=u+\nabla \pi_{h,B}$. Here,
\ben\label{eq:pressure decomposition}
\nabla \pi_{h,B}=-E_B(u),\quad \nabla \pi_{1,B}=-E_B(u\cdot\nabla u),\quad \nabla \pi_{2,B}=E_B(\Delta u).
\een
\end{Definition}
Noting the properties of the projection operator $E_B$ and using $\eqref{eq:bound of wolf}$, $\eqref{eq:bound of wolf'}$, there holds
\ben\label{ine:estimate p h}
\|\nabla \pi_{h,B}\|_{L^s(B)} \leq C \|u\|_{L^s(B)} \quad {\rm for} \quad s > 1,
\een
\ben\label{ine:estimate p 1}
\|\pi_{1,B}\|_{L^{s'}(B)} \leq C \|u\|_{L^{2s'}(B)}^2, \quad {\rm for} \quad s' > 1,
\een
and
\ben\label{ine:estimate p 2}
\|\pi_{2,B}\|_{L^2(B)} \leq C \|\nabla u\|_{L^2(B)}.
\een

\begin{Remark}
 Wolf  \cite{[Wolf10]} proved the existence of a local suitable weak solution to the Navier-Stokes equations (\ref{NSE}).
Chae and Wolf  \cite{CW2017} proved that if $(u,\pi)$ is a suitable weak solution to the Navier-Stokes equations (\ref{NSE}), then $u$ is a local suitable weak solution in the sense of Definition \ref{def:1}.
\end{Remark}

We also need the following Riesz potential estimate (see, for example, \cite[ p159]{GT}).
\begin{Lemma}\label{lem:Riesz}
	Let $\Omega$ be a bounded domain, $\mu\in(0,1]$, $1\leq q\leq \infty$, $0\leq \delta=1/p-1/q<\mu$ and $V_\mu f(x)=\int_{\Omega}|x-y|^{n(\mu-1)}f(y)dy$, then we have
\beno
\|V_\mu f(x)\|_{L^q(\Omega)}\leq \big(\frac{1-\delta}{\mu-\delta}\big)^{1-\delta}w_n^{1-\mu}|\Omega|^{\mu-\delta}\|f\|_{L^p(\Omega)},
\eeno
where $w_n$ is the volume of unit ball in $\mathbb{R}^n$.
\end{Lemma}

\section{Cacciopolli's inequality}
In this section, we establish a Cacciopolli's inequality for the Navier-Stokes equations in term of velocity only for the value range of p in $[2,\infty]$.
\begin{Proposition}\label{cacci}
Assume that  $(u,\pi)$ is a suitable weak solution to (\ref{NSE}) in $Q_1$.
	For any $(p,q)$ satisfying $1\leq 2/p+3/q<2$ with  $2\leq p\leq\infty$, $u\in L_t^pL_x^q(Q_1)$. Then the following Cacciopolli's inequality holds true:
\beno
\|u\|_{L^pL^{\frac{6p}{3p-4}}(Q_\frac34)}^2+\|\nabla u\|^2_{L^2L^2(Q_\frac34)}\leq C \|u\|_{L^pL^q(Q_1)}^2 + C \|u\|_{L^pL^q(Q_1)}^4 + C \|u\|_{L^pL^q(Q_1)}^{\frac{2\alpha}{\alpha-1}}.
\eeno
where $\alpha = \frac 2{\frac2p+\frac3q}$.
\end{Proposition}

%


\begin{Remark}
The case $2 \leq p < 3$ is specially, since $\|u\|_{L^3L^3}$ can not be controlled by $\|u\|_{L^p L^q}$ and energy norm for $p \in [2,3)$ with respect to time direction $t$.
It is worth mentioning that the energy norms of $v_B = u + \nabla \pi_h$ include all indicators for $p \geq 2$ with respect to time direction $t$ but it is not involved  into the iterative process since $\nabla \pi_h$ is related with the domain $B_\rho$.
\end{Remark}


\subsection{Proof of Proposition \ref{cacci}:  The case of $p > 3$}

Firstly, for $\frac34 \leq \varrho <\rho \le 1$, let  $Q_\rho = (-\rho^2,0) \times B_\rho$, and $B_\rho = \{x\in \mathbb{R}^3;|x|\leq \rho\}$.
Write $B = B_\rho$ and define $v_B = u + \nabla \pi_{h,B}$ with $\nabla \pi_{h,B} = -E_{B}(u)$.
Choose a cut-off function as
\beno
\phi = 1 \quad {\rm in} \quad Q_{\sigma_1} \quad {\rm with} \quad \sigma_1 = \frac{2\varrho+\rho}3,
\eeno
\beno
\phi = 1 \quad {\rm on} \quad Q_{\sigma_2}^c \quad {\rm with} \quad \sigma_2 = \frac{\varrho+2\rho}3,
\eeno
and satisfies
\ben\label{ine:phi nabla}
|\nabla \phi| \leq C (\rho-\varrho)^{-1}, \quad |\partial_t \phi| + |\nabla^2 \phi| \leq C (\rho-\varrho)^{-2}.
\een
Secondly, choosing $\varphi = \phi^2$ in the local energy inequality $\eqref{eq:local energy inq}$, we have
\ben\label{reloce} \nonumber
&&\int_{B_\rho}|v_B(x,s)|^2 \phi^2 dx+\int_{Q_\rho}|\nabla v_B(x,\tau)|^2 \phi^2 dxd\tau\\ \nonumber
&\leq& C(\rho-r)^{-2}\int_{Q_\rho}|v_B(x,\tau)|^2dxd\tau+C(\rho-r)^{-1}\int_{Q_\rho}|v_B|^2|v_B-\nabla \pi_{h,B}| dxd\tau\\ \nonumber
&&+C\int_{Q_{\sigma_2}} |v_B|^2 |\nabla \nabla \pi_{h,B}| dxd\tau+C\int_{Q_{\sigma_2}} |\nabla \pi_{h,B}| |\nabla \nabla \pi_{h,B}| |v_B| dxd\tau\\ \nonumber
&&+C(\rho-r)^{-1}\int_{Q_\rho} |\pi_{1,B}+\pi_{2,B}| |v_B|dxd\tau \\
&:=& I_1 + I_2 + I_3 + I_4 + I_5,
\een
where
\beno
\nabla \pi_{1,B} = -E_{B}(u\cdot\nabla u),\quad {\rm and} \quad \nabla \pi_{2,B}=E_{B}(\triangle u).
\eeno

{\bf Step I: \underline{Estimate of local energy inequality via $\|u\|_{L^3(Q_\rho)}$}.}
By (\ref{ine:estimate p h}), there holds for any $s \in (1,6)$,
\ben\label{ine:v_B u}
\|v_B\|_{L^s(B_\rho)} = \|u + \nabla \pi_{h,B}\|_{L^s(B_\rho)} \leq C \|u\|_{L^s(B_\rho)},
\een
which yields that
\beno
I_1 \leq C (\rho - \varrho)^{-2} \int_{Q_\rho} |u|^2 dx d\tau.
\eeno
Similarly, by $\eqref{ine:v_B u}$ and (\ref{ine:estimate p h}), the estimate of $I_2$ is that
\beno
I_2 \leq C (\rho-\varrho)^{-1} \int_{Q_\rho}|u|^3 dx d\tau.
\eeno

For the part of $\nabla \nabla \pi_{h,B}$, noting that $-\Delta v_h + \nabla \pi_h = -u$, which implies that $\Delta \pi_h = 0$, and $ \pi_h$ is harmonic.
Recall the estimates of harmonic function (see, for example, \cite{JWZ}): for any $1\leq p,q\leq \infty$, $0<\varrho<\rho$, and any harmonic function $h$, it holds
\ben\label{prop:harmonic}
\|\nabla^k h\|_{L^q(B_\varrho)}\leq \frac{C\varrho^{\frac3 q}}{(\rho-\varrho)^{\frac3p+k}}\|h\|_{L^p(B_\rho)}.
\een
Using $\eqref{ine:v_B u}$, (\ref{ine:estimate p h}), $\eqref{prop:harmonic}$ and the H\"{o}lder's inequality, for the term $I_3$, there holds
\beno
I_3 &=& C\int_{Q_{\sigma_2}} |v_B|^2 |\nabla \nabla \pi_{h,B}| dxd\tau \\
&\leq& C \left(\int_{Q_\rho} |v_B|^3 dxd\tau\right)^\frac23 \left(\int_{Q_{\sigma_2}} |\nabla \nabla \pi_{h,B}|^3 dxd\tau\right)^\frac13 \\
&\leq& C \left(\int_{Q_\rho} |u|^3 dxd\tau\right)^\frac23 \left(C\sigma_2^3 (\rho-\varrho)^{-6}\int_{Q_\rho} |\nabla \pi_{h,B}|^3 dxd\tau\right)^\frac13\\
&\leq& C \rho (\rho-\varrho)^{-2} \int_{Q_\rho}|u|^3 dx d\tau.
\eeno
Similarly,
\beno
I_4 = C\int_{Q_{\sigma_2}} |v_B| |\nabla \pi_h| |\nabla \nabla \pi_h| dxd\tau \leq C \rho (\rho-r)^{-2} \int_{Q_\rho}|u|^3 dx d\tau.
\eeno

For the term $I_5$, using $\eqref{ine:v_B u}$, (\ref{ine:estimate p 1}) and (\ref{ine:estimate p 2}), we infer that
\beno
I_5 &=& C(\rho-\varrho)^{-1}\int_{Q_\rho} |\pi_{1,B}+\pi_{2,B}| |v_B|dxd\tau \\
&\leq& C (\rho-\varrho)^{-1} \left(\int_{Q_\rho} |v_B|^3 dxd\tau\right)^\frac13 \left(\int_{Q_\rho} |\pi_{1,B}|^\frac32dxd\tau\right)^\frac23 \\
&& \quad \quad +  C (\rho-\varrho)^{-1} \left(\int_{Q_\rho} |v_B|^2 dxd\tau\right)^\frac12 \left(\int_{Q_\rho} |\pi_{2,B}|^2dxd\tau\right)^\frac12 \\
&\leq& C (\rho-\varrho)^{-1} \int_{Q_\rho} |u|^3 dx d\tau + C (\rho-\varrho)^{-1} \left(\int_{Q_\rho} |u|^2 dx d\tau\right)^\frac12 \left(\int_{Q_\rho}|\nabla u|^2 dx d\tau\right)^\frac12 \\
&\leq& C (\rho-\varrho)^{-1} \int_{Q_\rho} |u|^3 dx d\tau + C (\rho-\varrho)^{-2} \int_{Q_\rho} |u|^2 dx d\tau + \frac1{16} \int_{Q_\rho}|\nabla u|^2 dx d\tau.
\eeno

Collecting the term $I_1 - I_5$, we arrive at
\ben\label{ine:energy}\nonumber
&&\int_{B_\rho}|v_B(x,s)|^2 \phi dx+\int_{Q_\rho}|\nabla v_B(x,\tau)|^2 \phi dxd\tau\\ \nonumber
&\leq& C(\rho-\varrho)^{-2}\int_{Q_\rho}|u|^2+C\left((\rho-\varrho)^{-1}+\rho(\rho-\varrho)^{-2}\right)\int_{Q_\rho}|u|^3+\frac{1}{16}\|\nabla u\|^2_{L^2L^2(Q_\rho)}\\
&:= &J_1 + J_2 + J_3.
\een

{\bf Step II: \underline{Estimate the term $\|u\|_{L^3(Q_\rho)}$}.}
It follows from  H\"{o}lder's inequality  that
\begin{align*}
J_1 \leq C \rho^{\frac53} (\rho-\varrho)^{-2} \|u\|_{L^3(Q_\rho)}^2 \leq C (\rho-\varrho)^{-2} \|u\|_{L^3(Q_\rho)}^2,
\end{align*}
since $\rho \leq 1$.

Next we deal with the term $\|u\|_{L^3(Q_\rho)}$.
For any $1\leq\tau\leq\frac p2$ and $\frac{3p}{3p-4}\leq\kappa\leq3$ satisfying $\frac{2}{\tau}+\frac{3}{\kappa}=3$, by interpolation inequality, we see that for any $f$,
\ben\label{simpleinterplation}
\|f^2\|_{L^{\tau}L^{\kappa}(Q_\rho)}=\|f\|^2_{L^{2\tau}L^{2\kappa}(Q_\rho)}
\leq C \left(\|f\|^2_{L^{p}L^{\frac{6p}{3p-4}} (Q_\rho)}+\|f\|_{L^{2}L^{6}(Q_\rho)}^2\right).
\een
For  a fixed $p>3$, let $q\in \left(\frac{3}{2-\frac{2}{p}}, \frac{3}{\frac{2}{p-2}-\frac{2}{p}}\right]$ and $\alpha (\frac{2}{p}+\frac{3}{q})=2$, then $1<\alpha \leq p-2$. For $1<\alpha \leq \min\{2,p-2\}$, let
\beno
 r = \frac p\alpha = \frac p 2(\frac2p+\frac3q), \quad s = \frac q \alpha = \frac q2(\frac2p+\frac3q),
\eeno
then
\beno
\frac 1s + \frac 1{s'} = 1, \quad \frac 1r + \frac 1{r'} = 1, \quad \frac 2{r'} + \frac 3{s'} = 3,
\eeno
where, $s',r'$ are the conjugate index of $s,r$. Noting that
\beno
\iint_{Q_\rho}|u|^{3}dxdt=\iint_{Q_\rho}|u|^{\alpha}|u|^{3-\alpha}dxdt
				\leq\||u|^{\alpha}\|_{L^rL^s(Q_\rho)}\||u|^{3-\alpha}\|_{L^{r'}L^{s'}(Q_\rho)},
\eeno
by \eqref{simpleinterplation}, we get
\beno
\int_{Q_\rho}|u|^{3}dxdt &\leq& C \|u\|_{L^{\alpha r}L^{\alpha s}(Q_\rho)}^{\alpha}\Big(\||u|^{\frac{3-\alpha}2}\|_{L^pL^{\frac{6p}{3p-4}}(Q_\rho)}^{2}+\||u|^{\frac{3-\alpha}2}\|_{L^2L^6(Q_\rho)}^{2}\Big)\\
&\leq& C \|u\|_{L^{\alpha r}L^{\alpha s}(Q_\rho)}^{\alpha}\Big(\|u\|^{3-\alpha}_{L^{\frac{3-\alpha}2 p} L^{\frac{3-\alpha}2\frac{6p}{3p-4}}(Q_\rho)}+\|u\|^{3-\alpha}_{L^{3-\alpha}L^{9-3\alpha}(Q_\rho)}\Big),
\eeno
since $\frac p{p-2}\leq r\leq \infty$ and $r'\in [1,p/2]$ due to $\alpha\leq p-2$. Then H\"{o}lder's inequality implies that
\beno
\|u\|^{3-\alpha}_{L^{\frac{3-\alpha}2 p} L^{\frac{3-\alpha}2 \frac{6p}{3p-4}}(Q_\rho)}+\|u\|^{3-\alpha}_{L^{3-\alpha}L^{9-3\alpha}(Q_\rho)} \leq C \rho^{\frac{3}{2}(\alpha-1)}\left(  \|u\|_{L^{p}L^{\frac{6p}{3p-4}}(Q_\rho)}^{3-\alpha}+\|\nabla u\|_{L^{2}(Q_\rho)}^{3-\alpha}\right),
\eeno
which means
\ben\label{ine:u 3 energy}
\int_{Q_\rho}|u|^{3}dxdt \leq C\|u\|_{L^pL^q(Q_\rho)}^{\alpha}\Big(\|u\|_{L^{p}L^{\frac{6p}{3p-4}}(Q_\rho)}^{2}+\|\nabla u\|_{L^{2}(Q_\rho)}^{2}\Big)^{\frac{3-\alpha}2}.
\een
In addition, for the case of $q>\frac{3}{\frac{2}{p-2}-\frac{2}{p}}$, the estimate of (\ref{ine:u 3 energy}) still holds due to the   H\"{o}lder's inequality.

Using $\eqref{ine:u 3 energy}$ and Young's inequality, there holds
\ben\label{ine:J1}
J_1 \leq \frac1{32} \left(\|u\|_{L^pL^{\frac{6p}{3p-4}}(Q_\rho)}^2 + \|\nabla u\|_{L^2(Q_\rho)}^2\right) + C (\rho-\varrho)^{-\frac 6\alpha} \|u\|_{L^pL^q(Q_\rho)}^2.
\een
Similarly,
\ben\label{ine:J2}
J_2 \leq \frac1{32} \left(\|u\|_{L^pL^{\frac{6p}{3p-4}}(Q_\rho)}^2 + \|\nabla u\|_{L^2(Q_\rho)}^2\right) + C (\rho-\varrho)^{-\frac{4}{\alpha-1}}  \|u\|_{L^pL^q(Q_\rho)}^{\frac{2\alpha}{\alpha-1}}.
\een
Combining $\eqref{ine:energy}$, $\eqref{ine:J1}$ and $\eqref{ine:J2}$, we arrive at
\ben\label{ine:energy estimate}\nonumber
&&\int_{B_\rho} |v_B|^2 \phi^2 + \int_{Q_\rho} |\nabla v_B|^2 \phi^2 \leq \frac18 \left(\|u\|_{L^pL^{\frac{6p}{3p-4}}(Q_\rho)}^2 + \|\nabla u\|_{L^2(Q_\rho)}^2\right) \\
&&+ C (\rho-\varrho)^{-\frac 6\alpha} \|u\|_{L^pL^q(Q_\rho)}^2 + C (\rho-\varrho)^{-\frac{4}{\alpha-1}}  \|u\|_{L^pL^q(Q_\rho)}^{\frac{2\alpha}{\alpha-1}}.
\een

{\bf Step III:  \underline{Estimate the terms including $v_B$}.}
Noting that $\frac2p+\frac3{\frac{6p}{3p-4}} = \frac32$, using Young's inequality and Sobolev's embedding, there holds
\ben\label{ine:energy v B phi}
\|v_B\|_{L^pL^{\frac{6p}{3p-4}}(Q_\varrho)}^2 \leq C \int_{B_\varrho} |v_B(t)|^2 dx + C \int_{Q_\varrho} |\nabla v_B|^2 dx dt.
\een
Since $\nabla \pi_{h,B}$ is harmonic function, using (\ref{prop:harmonic}), there holds
\ben\label{ine:energy p h phi} \nonumber
\|\nabla \pi_h \|_{L^pL^{\frac{6p}{3p-4}}(Q_{\sigma_1})}^2 &\leq& C \rho^{3-\frac4p} (\rho-\varrho)^{-\frac6q} ||\nabla \pi_h||_{L^pL^q(Q_{\sigma_2})}^2 \\
&\leq& C \rho^{3-\frac4p} (\rho-\varrho)^{-\frac6q} ||u||_{L^pL^q(Q_{\rho})}^2.
\een
Using (\ref{ine:energy estimate}), (\ref{ine:energy v B phi}) and (\ref{ine:energy p h phi}), we conclude by the triangle inequality that
\ben\label{ine:u p 6p} \nonumber
\|u\|_{L^pL^{\frac{6p}{3p-4}}(Q_\varrho)}^2 &\leq& 2 \|v_B\|_{L^pL^{\frac{6p}{3p-4}}(Q_\varrho)}^2 + 2\|\nabla \pi_{h,B}\|_{L^pL^{\frac{6p}{3p-4}}(Q_{\varrho})}^2\\ \nonumber
&\leq& C\|v_B \phi\|_{L^\infty L^{2}(Q_\rho)}^2 + C \|\nabla v_B \phi\|_{L^2(Q_\rho)}^2 + 2\|\nabla \pi_{h,B}\|_{L^pL^{\frac{6p}{3p-4}}(Q_{\sigma_1})}^2\nonumber \\ \nonumber
&\leq& \frac 14 \left(\|u\|_{L^pL^{\frac{6p}{3p-4}}(Q_\rho)}^2 + \|\nabla u\|_{L^2(Q_\rho)}^2\right) + C (\rho-\varrho)^{-\frac 6\alpha} \|u\|_{L^pL^q(Q_\rho)}^2 \\ \nonumber
&&+ C \left((\rho-\varrho)^{-1}+(\rho-\varrho)^{-2}\right)^\frac{2}{\alpha-1} \|u\|_{L^pL^q(Q_\rho)}^{\frac{2\alpha}{\alpha-1}} \\
&&+ C (\rho-\varrho)^{-\frac6q} ||u||_{L^pL^q(Q_\rho)}^2.
\een
Similarly, noting that $\nabla u = \nabla v_B - \nabla \nabla \pi_{h,B}$, for almost $t \in I_\rho$, $\eqref{prop:harmonic}$ and $\eqref{ine:v_B u}$ implies
\beno
&&\|\phi(t) \nabla u(t)\|_{L^2(B_\rho)}^2 \\
&=& \int_{B_\rho} |\phi(t) \nabla v_B(t)|^2 dx - \int_{B_\rho} (\nabla v_B(t) + \nabla u(t)) : (\nabla v_B(t) - \nabla u(t)) \phi^2(t) dx \\
&=& \int_{B_\rho} |\phi(t) \nabla v_B(t)|^2 dx - \int_{B_\rho} (\nabla v_B(t) + \nabla u(t)) : \nabla \nabla \pi_{h,B} \phi^2(t) dx \\
&=& \int_{B_\rho} |\phi(t) \nabla v_B(t)|^2 dx + \int_{B_\rho} (v_B(t) + u(t)) \cdot \nabla \nabla \pi_{h,B} \cdot \nabla \phi^2(t) dx \\
&\leq& \int_{B_\rho} |\phi(t) \nabla v_B(t)|^2 dx + C (\rho - r)^{-1} \left(\int_{B_\rho} |v_B + u|^2 dx\right)^\frac12 \left(\int_{B_{\sigma_2}} |\nabla^2 \pi_{h}|^2 dx\right)^\frac12 \\
&\leq& \int_{B_\rho} |\phi(t) \nabla v_B(t)|^2 dx + C (\rho - r)^{-2} \int_{B_\rho} |u|^2 dx.
\eeno
Integrating with respect to $t$, we have
\ben\label{ine:nabla u 2}
\int_{Q_\rho} |\nabla u \phi|^2 \leq \int_{Q_\rho} |\phi \nabla v_B|^2 dx + C (\rho - r)^{-2} \int_{Q_\rho} |u|^2 dx.
\een

Then for $\frac34 \leq \varrho < \rho\leq 1$, $\alpha=\frac2{\frac2q + \frac3p}$, combining (\ref{ine:u p 6p}), (\ref{ine:nabla u 2}), (\ref{ine:energy estimate}) and $\eqref{ine:J1}$, we have
\beno
\|u\|_{L^pL^{\frac{6p}{3p-4}}(Q_\varrho)}^2 + \|\nabla u\|_{L^2(Q_\varrho)}^2 &\leq& \frac 34 (\|u\|_{L^pL^{\frac{6p}{3p-4}}(Q_\rho)}^2 + \|\nabla u\|_{L^2(Q_\rho)}^2) \\
&&+ C (\rho-\varrho)^{-\frac 6\alpha} \|u\|_{L^pL^q(Q_\rho)}^2 + C (\rho-\varrho)^{-\frac6q} ||u||_{L^pL^q(Q_\rho)}^2\\
&&+ C (\rho-\varrho)^{-\frac{4}{\alpha-1}} \|u\|_{L^pL^q(Q_\rho)}^{\frac{2\alpha}{\alpha-1}}.
\eeno
Applying the iteration lemma (see \cite[Lemma V.3.1,   p.161 ]{[Giaquinta]}), we end up with
\beno
\|u\|_{L^pL^{\frac{6p}{3p-4}}(Q_\varrho)}^2 + \|\nabla u\|_{L^2(Q_\varrho)}^2 &\leq&  C  (\rho-\varrho)^{-\frac 6\alpha} \|u\|_{L^pL^q(Q_\rho)}^2 + C (\rho-\varrho)^{-\frac6q} ||u||_{L^pL^q(Q_\rho)}^2\\
&&+  C (\rho-\varrho)^{-\frac{4}{\alpha-1}} \|u\|_{L^pL^q(Q_\rho)}^{\frac{2\alpha}{\alpha-1}},
\eeno
which leads to that
\beno
\|u\|_{L^pL^{\frac{6p}{3p-4}}(Q_\frac34)}^2 + \|\nabla u\|_{L^2(Q_\frac34)}^2 \leq C ||u||_{L^pL^q(Q_1)}^2 + C ||u||_{L^pL^q(Q_1)}^{\frac{2\alpha}{\alpha-1}}.
\eeno
The proof is complete.

\subsection{Proof of Proposition \ref{cacci}:  The case of  $2\leq  p \leq 3$}

Let $\frac34 \leq \varrho < \rho\leq 1$ and  we still use the cut-off function $\phi$ as in the last subsection.
Taking the test function $\varphi=\phi^{2\beta}$in the local energy inequality (\ref{eq:local energy inq}), we have
\beno
&&\|v_B\phi^\beta\|^2_{L^\infty L^2(Q_{\rho})}+2\|\nabla (v_B\phi^\beta)\|^2_{L^2 L^2(Q_\rho)}\\
&\leq& C \left|\int_{Q_\rho} |v_B|^2 (\partial_t(\phi^{2\beta})+\Delta(\phi^{2\beta})+(\nabla \phi)^2\phi^{2\beta-2})\right|  + \int_{Q_\rho}u\cdot\nabla (\phi^{2\beta}) |v_B|^2 \\
&& +2\int_{Q_\rho} u \cdot \nabla \nabla \pi_{h,B} \cdot v_B \phi^{2\beta} +2\int_{Q_\rho} \pi_{1,B}v_B\cdot \nabla(\phi^{2\beta}) +2\int_{Q_\rho} \pi_{2,B} v_B \cdot \nabla(\phi^{2\beta}) \\
&:=& K_1 + K_2 + K_3 + K_4 + K_5.
\eeno

{\bf \underline{Estimate of  $K_1$}.}
Noting that $q > \frac 94$ for $p \in [2,3]$ since $\frac 2p + \frac 3q < 2$, using $\eqref{ine:v_B u}$ and $\eqref{ine:phi nabla}$, there hold
\ben\label{eq:K1}
K_1 \leq C (\rho-\varrho)^{-2} \int_{Q_\rho} |u|^2 dx dt \leq C (\rho-\varrho)^{-2} \|u\|_{L^p L^q(Q_\rho)}^2.
\een

{\bf \underline{Estimate of  $K_2$}.}
Taking $\beta = \beta_0$ large enough such that $2\beta_0-1\geq (3-\alpha)\beta_0$ with $\alpha=\frac2{\frac2p+\frac3q}$, using H\"{o}lder's inequality, we have
\ben\label{eq:K2}
K_2& \leq& C (\rho-\varrho)^{-1} \iint_{Q_\rho} |u||v_B|^2\phi^{2\beta-1}dxdt\nonumber \\
&\leq& C (\rho-\varrho)^{-1} \|u\|_{L^pL^q(Q_\rho)}\||v_B|^{\alpha-1}\|_{L^{\frac p{\alpha-1}}L^{\frac q{\alpha-1}}(Q_\rho)}
\||v_B\phi^\beta|^{3-\alpha}\|_{L^{r}L^{s}(Q_\rho)}.
\een
Here,
\beno
r = \frac p{p-\alpha}, \quad s = \frac q{q-\alpha}, \quad \alpha = \frac 2{\frac2p+\frac3q}.
\eeno
Claim that:   $r(3-\alpha) \geq 2$ and $2 \leq s(3-\alpha) \leq 6$ for $q\in \left(\frac{3}{2-\frac{2}{p}},\frac{7}{2-\frac{2}{p}}\right]$.
First, since $p \leq 3$, we have
\beno
r(3-\alpha) \geq  3.
\eeno
Second, due to $q\leq \frac{7}{2-\frac{2}{p}}$ there holds
\beno
s(3-\alpha) =(3-\alpha)\frac q{q-\alpha} \geq 2,
\eeno
and
\beno
s(3-\alpha) \leq 6,
\eeno
since $q>\frac94.$

Noting that
\beno
\frac2{r(3-\alpha)}+\frac3{s(3-\alpha)} = \frac3{3-\alpha} > \frac32,
\eeno
with $\frac2r+\frac3s=3$,
by H\"{o}lder's inequality and $\eqref{ine:v_B u}$ again, there holds
\ben\label{eq:K2 estimate}
K_2 &\leq& C (\rho-r)^{-1} \rho^{\frac32(\alpha-1)} \|u\|_{L^pL^q}\|v_B\|_{L^{p}L^{q}(Q_\rho)}^{\alpha-1}
\Big(\|v_B\phi^\beta\|_{L^{\infty}L^{2}(Q_\rho)}^{2}+\|\nabla (v_B\phi^\beta)\|_{L^{2}(Q_\rho)}^{2}\Big)^{\frac{3-\alpha}2} \nonumber\\
& \leq &\frac 18 \Big(\|v_B\phi^\beta\|_{L^{\infty}L^{2}(Q_\rho)}^{2}+\|\nabla (v_B\phi^\beta)\|_{L^{2}(Q_\rho)}^{2}\Big) + C (\rho-\varrho)^{-\frac2{\alpha-1}} \|u\|_{L^pL^q(Q_\rho)}^{\frac{2\alpha}{\alpha-1}}.
\een
In addition, for the case of $q>\frac{7}{2-\frac{2}{p}}$, one can take $q_0=\frac{7}{2-\frac{2}{p}}$ and (\ref{eq:K2 estimate}) holds for $q_0$. Then the H\"{o}lder inequality applies and (\ref{eq:K2 estimate}) holds for and $q>q_0$.


{\bf \underline{Estimate of  $K_3$}.} Noting that $q > \frac94$ since $p \in [2,3]$ and $\rho\leq1$, using (\ref{prop:harmonic}), $\eqref{ine:v_B u}$ and H\"{o}lder's inequality, we know that
\ben\label{eq:K3}
K_3 &\leq& C \rho^{\frac72 - \frac 3q - \frac 4p} \|v_B \phi^\beta\|_{L^\infty L^2(Q_\rho)} \|u\|_{L^p L^q(Q_\rho)} \|\nabla \nabla \pi_{h,B}\|_{L^p L^\infty(Q_{\sigma_2})}\nonumber \\
& \leq& C (\rho-\varrho)^{-1-\frac3q} \|v_B \phi^\beta\|_{L^\infty L^2(Q_\rho)} \|u\|_{L^p L^q(Q_\rho)} \|u\|_{L^p L^q(Q_{\rho})}\nonumber \\
& \leq& \frac 18 \Big(\|v_B\phi^\beta\|_{L^{\infty}L^{2}(Q_\rho)}^{2}+\|\nabla (v_B\phi^\beta)\|_{L^{2}(Q_\rho)}^{2}\Big) \nonumber \\
&&+ C (\rho-\varrho)^{-2-\frac6q} \|u\|_{L^pL^q(Q_\rho)}^4.
\een

{\bf \underline{Estimate of  $K_4$}.} Recall that
\beno
K_4 = -4 \beta \int_{Q_\rho} v_B \phi^\beta \cdot \nabla \phi ~ \pi_{1,B} \phi^{\beta-1} dxds.
\eeno
Now we rewrite the first equation of (\ref{NSE}) as
\beno
\partial_tv_B-\Delta v_B+u\cdot\nabla u+\nabla \pi_{1,B}+\nabla \pi_{2,B}=0,
\eeno
with $v_B=u+\nabla \pi_{h,B}$ and $\nabla \pi_{h,B} = -E_{B_\rho}(u)$. Using the representation formula of pressure $\pi_{1,B}$, we have
\ben\label{localpre}
    \pi_{1,B}\xi&=&R_iR_j(\xi u_iu_j)-N*(\partial_{ij}\xi u_iu_j)+\partial_jN*(u_iu_j\partial_i\xi)\nonumber \\
&&+\partial_iN*(u_iu_j\partial_j\xi)-N*(\pi_{1,B}\Delta \xi)+2\partial_jN*(\partial_j\xi \pi_{1,B}),
\een
where $\xi$ is a cutoff function,  $N=-\frac 1{4\pi|x|}$ is the kernel of Poisson equation and $R_i = \frac{\partial_i}{\sqrt{-\Delta}}$ is Riesz transform. Rewrite $K_4 = -4 \beta (K_{41} + \cdots + K_{46}).$

Choose $\xi=\phi^{\beta-1}$,  and note that $\beta =\beta_0$ satisfying  $\beta_0-1\geq (2-\alpha)\beta_0$. Then
\beno
K_{41} &=& \int_{Q_\rho} v_B\phi^\beta \cdot \nabla \phi R_iR_j(\xi u_iu_j) dxdt\nonumber\\
&=&\int\big( R_iR_j(\xi u_i(v_B)_j)-R_iR_j(\xi u_i(\nabla \pi_{h,B})_j)\big) \phi^\beta v_B\cdot \nabla \phi dxdt
\eeno
Using Lemma \ref{lem:Riesz}, the same estimates as $K_2$ in (\ref{eq:K2}) and $K_3$ in (\ref{eq:K3}) yields that
\ben\label{eq:K41}
K_{41}&\leq&
C (\rho - r)^{-1} \|v_B\phi^\beta\|_{L^{(3-\alpha)r}L^{(3-\alpha)s}(Q_\rho)}\|u\|_{L^pL^q(Q_\rho)}
\|v_B\|^{\alpha-1}_{L^pL^q(Q_\rho)}\|v_B\phi^\beta\|^{2-\alpha}_{L^{(3-\alpha)r}L^{(3-\alpha)s}(Q_\rho)}\nonumber\\
& &\quad +C (\rho - r)^{-1} \|v_B\phi^\beta\|_{L^\infty L^2(Q_\rho)}\|u\|_{L^pL^q(Q_\rho)}\|\nabla \pi_{h,B}\|_{L^pL^\infty(Q_{\sigma_2})} \nonumber\\
&\leq& \frac 18 \Big(\|v_B\phi^\beta\|_{L^{\infty}L^{2}(Q_\rho)}^{2}+\|\nabla (v_B\phi^\beta)\|_{L^{2}(Q_\rho)}^{2}\Big) + C (\rho-r)^{-\frac2{\alpha-1}} \rho^3 \|u\|_{L^pL^q(Q_\rho)}^{\frac{2\alpha}{\alpha-1}} \nonumber\\
&& \quad + C (\rho-r)^{-2-\frac6q} \|u\|_{L^pL^q(Q_\rho)}^4.
\een

For $K_{42}$, using H\"{o}lder's inequality and Young's inequality for convolution form, there holds
\beno
K_{42} &\leq& C (\rho-r)^{-1} \|v_B \phi^\beta\|_{L^\infty L^2 (Q_\rho)} \|N*(\partial_{ij}\xi u_iu_j)\|_{L^1 L^2(Q_{\sigma_2})} \\
&\leq& C (\rho-r)^{-1} \|v_B \phi^\beta\|_{L^\infty L^2 (Q_\rho)} (\rho-r)^{-2} \|N\|_{L^\infty L^2(Q_\ast)} \||u|^2\|_{L^1L^1(Q_{\sigma_2})}
\eeno
where   $Q_\ast = (-\rho^2,0) \times B_\ast$, $B_\ast = \{x:|x|\leq 2\sigma_2\}$. It follows that
\ben\label{ine:N}
\|N\|_{L^\infty L^2(Q_\ast)} \leq C \rho^\frac12\leq C.
\een
Noting that $q\in(\frac94,+\infty)$ since $p\in[2,3]$, 
using $\eqref{ine:N}$, we deduce
\ben\label{eq:K42} \nonumber
K_{42} &\leq& C \rho^\frac12 (\rho-\varrho)^{-3} \|v_B \phi^\beta\|_{L^\infty L^2 (Q_\rho)} \|u\|_{L^2 L^2(Q_\rho)}^2 \\
&\leq& C (\rho-\varrho)^{-3} \|v_B \phi^\beta\|_{L^\infty L^2 (Q_\rho)} \|u\|_{L^p L^q(Q_\rho)}^2.
\een

For the term of $K_{43}$, we have
\ben\label{eq:K43}
K_{43}\leq C (\rho-\varrho)^{-2} \left(\|v_B \phi^\beta\|_{L^\infty L^2 (Q_\rho)}+\|\nabla (v_B\phi^\beta)\|_{L^2 L^2(Q_\rho)}\right) \|u\|_{L^p L^q(Q_\rho)}^2.
\een
 First, when $p = 2$, i follows that  $q > 3$ since $\frac 2p+\frac 3q < 2$.
Using H\"{o}lder's inequality and Young's inequality, we arrive at
\beno
K_{43} &\leq &C (\rho-\varrho)^{-1} \|v_B \phi^\beta\|_{L^\infty L^2(Q_\rho)} \|\partial_jN*(u_iu_j\partial_i\xi)\|_{L^1 L^2(Q_{\sigma_2})} \\
&\leq& C (\rho-\varrho)^{-2} \|v_B \phi^\beta\|_{L^\infty L^2(Q_\rho)} \|\partial_j N\|_{L^\infty L^\frac{2q}{3q-4}(Q_\ast)} \||u|^2\|_{L^1 L^\frac q2(Q_\rho)}\\
&\leq& C (\rho-\varrho)^{-2} \left(\|v_B \phi^\beta\|_{L^\infty L^2 (Q_\rho)}+\|\nabla (v_B\phi^\beta)\|_{L^2 L^2(Q_\rho)}\right) \|u\|_{L^p L^q(Q_\rho)}^2,
\eeno
where  the term $\|\partial_j N\|_{L^\infty_t L^\frac{2q}{3q-4}_x(Q_\rho)}$ is integrable since $\frac{2q}{3q-4} < \frac32$ when $q\in (3,4]$. When $q>4$, the above inequality holds for $q_0=4$, then it is still true for $q>4$ by the H\"{o}lder inequality. Second,
For the case $2<p\leq3$,  using  H\"{o}lder's inequality and Young's inequality again, we get
\beno
K_{43} &\leq& C (\rho-\varrho)^{-1} \|v_B \phi^\beta\|_{L^\frac p{p-2} L^\frac{6p}{8-p} (Q_\rho)} \|\partial_jN*(u_iu_j\partial_i\xi)\|_{L^\frac p2 L^\frac{6p}{7p-8}(Q_{\sigma_2})} \\
&\leq& C (\rho-\varrho)^{-2} \|v_B \phi^\beta\|_{L^\frac p{p-2} L^\frac{6p}{8-p} (Q_\rho)} \|\partial_j N\|_{L^\infty L^\frac 1{\frac{13}6 - \frac23(\frac2p+\frac3q)}(Q_\ast)} \||u|^2\|_{L^\frac p2 L^{\frac{q}{2}}(Q_\rho)}.
\eeno
Obviously, the term $\|\partial_j N\|_{L^\infty L^\frac 1{\frac{13}6 - \frac23(\frac2p+\frac3q)}(Q_\ast)}$ is integrable since
\beno
1\leq \frac 1{\frac{13}6 - \frac23(\frac2p+\frac3q)} < \frac32, \quad q\leq \frac{3}{\frac74-\frac2p}
\eeno
Note that
$\frac{2(p-2)} p + \frac{8-p}{2p} = \frac32$, where  $\frac p{p-2} \geq 2$ and $2\le \frac{6p}{8-p}\le 6$. Thus,
\beno
 \|v_B \phi^\beta\|_{L^\frac p{p-2} L^\frac{6p}{8-p} (Q_\rho)}\leq C\left(\|v_B \phi^\beta\|_{L^\infty L^2 (B_\rho)}+\|\nabla (v_B\phi^\beta)\|_{L^2 L^2(B_\rho)}\right)
\eeno
Thus (\ref{eq:K43}) is proved for $q\leq\frac{3}{\frac74-\frac2p}$. When $q>\frac{3}{\frac74-\frac2p}$, the above inequality holds for $q_0=\frac{3}{\frac74-\frac2p}$, then it is still true for $q>\frac{3}{\frac74-\frac2p}$ by the H\"{o}lder inequality. The proof of (\ref{eq:K43}) is complete.

Similarly,
\ben\label{eq:K44}
K_{44} \leq C (\rho-\varrho)^{-2} \left(\|v_B \phi^\beta\|_{L^\infty L^2 (B_\rho)}+\|\nabla (v_B\phi^\beta)\|_{L^2 L^2(B_\rho)}\right) \|u\|_{L^p L^q(Q_\rho)}^2.
\een

For $K_{45}$, noting that $q>\frac94$ and using (\ref{ine:N}), there holds
\ben\label{eq:K45} \nonumber
K_{45} &\leq& C (\rho-\varrho)^{-1} \|v_B \phi^\beta\|_{L^\infty L^2 (Q_\rho)} \|N*(\pi_{1,B}\Delta \xi)\|_{L^1 L^2(Q_{\rho})} \\ \nonumber
&\leq& C (\rho-\varrho)^{-1} \|v_B \phi^\beta\|_{L^\infty L^2 (Q_\rho)} (\rho-r)^{-2} \|N\|_{L^\infty L^2(Q_\ast)} \|\pi_{1,B}\|_{L^1L^1(Q_\rho)} \\ \nonumber
&\leq& C \rho^{\frac 12} \rho^{5-\frac6q-\frac4p} (\rho-\varrho)^{-3} \|v_B \phi^\beta\|_{L^\infty L^2 (Q_\rho)} \|\pi_{1,B}\|_{L^\frac p2 L^\frac q2(Q_\rho)} \\
&\leq& C (\rho-\varrho)^{-3} \|v_B \phi^\beta\|_{L^\infty L^2 (Q_\rho)} \|u\|_{L^p L^q(Q_\rho)}^2.
\een

For $K_{46}$, H\"{o}lder's inequality and Young's inequality, we deduce that
\beno
K_{46} &\leq& C (\rho-\varrho)^{-1}  \|v_B \phi^\beta\|_{L^\frac p{p-2} L^\frac{6p}{8-p} (Q_\rho)} \|\partial_jN*(\pi_{1,B} \partial_j\xi)\|_{L^\frac p2 L^\frac{6p}{7p-8}(Q_{\rho})}\\
&\leq&  C (\rho-\varrho)^{-1}  \|v_B \phi^\beta\|_{L^\frac p{p-2} L^\frac{6p}{8-p} (Q_\rho)} \|\partial_jN\|_{L^\infty L^\frac 1{\frac{13}6 - \frac23(\frac2p+\frac3q)}(Q_\ast)} \|\pi_{1,B} \partial_j\xi\|_{L^\frac p2 L^{\frac{q}{2}}(Q_{\rho})},
\eeno
which is similar as $K43$ and $\frac p{p-2}=\infty$ for $p=2.$ The same arguments yields that
\ben\label{eq:K46}
K_{46}\leq C (\rho-\varrho)^{-2} \left(\|v_B \phi^\beta\|_{L^\infty L^2 (Q_\rho)}+\|\nabla (v_B \phi^\beta)\|_{L^2 L^2(Q_\rho)}\right) \|u\|_{L^p L^{q}(Q_\rho)}^2.
\een

Combining $\eqref{eq:K41}$, $\eqref{eq:K42}$, $\eqref{eq:K43}$, $\eqref{eq:K44}$, $\eqref{eq:K45}$ and $\eqref{eq:K46}$, choosing $\beta =\beta_0$ and using Young's inequality, we have
\ben\label{ine:K4} \nonumber
&K_4 \leq \frac 1{16} \Big(\|v_B\phi^\beta\|_{L^{\infty}L^{2}(Q_\rho)}^{2}+\|\nabla (v_B\phi^\beta)\|_{L^{2}(Q_\rho)}^{2}\Big) + C (\rho-\varrho)^{-\frac2{\alpha-1}} \|u\|_{L^pL^q(Q_\rho)}^{\frac{2\alpha}{\alpha-1}} \\
& \quad +  C (\rho - \varrho)^{-6} \|u\|_{L^pL^q(Q_\rho)}^4.
\een

{\bf \underline{Estimate of  $K_5$}.} Using $\eqref{ine:estimate p 2}$ and integration by parts, there holds
\ben\label{ine:K5} \nonumber
&K_5 \leq C \rho^{\frac 52 - \frac 2p - \frac3q} (\rho-\varrho)^{-1}\|v_B\|_{L^pL^q(Q_\rho)} \|\pi_{2,B}\|_{L^2L^2(Q_\rho)} \\
&\quad \leq C (\rho-\varrho)^{-2}\|u\|_{L^pL^q(Q_\rho)}^2 + \frac14 \|\nabla u\|_{L^2L^2(Q_\rho)}^2.
\een

Combining $\eqref{eq:K1}$, $\eqref{eq:K2 estimate}$, $\eqref{eq:K3}$, $\eqref{ine:K4}$ and $\eqref{ine:K5}$, we have
\ben\label{ine:energy''} \nonumber
&&\|v_B\phi^\beta\|^2_{L^\infty L^2(Q_\rho)}+\|\nabla (v_B\phi^\beta)\|^2_{L^2 L^2(Q_\rho)} \leq C (\rho-r)^{-\frac2{\alpha-1}} \|u\|_{L^pL^q(Q_\rho)}^{\frac{2\alpha}{\alpha-1}} \\
&& \quad +  C (\rho - r)^{-8} \|u\|_{L^pL^q(Q_\rho)}^4 + C (\rho-r)^{-2}\|u\|_{L^pL^q(Q_\rho)}^2 + \frac14 \|\nabla u\|_{L^2L^2(Q_\rho)}^2.
\een
Noting that $\frac2p+\frac3{\frac{6p}{3p-4}} = \frac32$ with $p\geq 2$  and $q>\frac94$, we have
\ben\label{ine:energy'}\nonumber
&&\|v_B \phi^\beta\|_{L^p L^\frac{6p}{3p-4}(Q_\rho)} + \|\nabla (v_B\phi^\beta)\|^2_{L^2 L^2(Q_\rho)} \leq C (\rho-r)^{-\frac2{\alpha-1}} \|u\|_{L^pL^q(Q_\rho)}^{\frac{2\alpha}{\alpha-1}} \\
&& \quad +  C (\rho - r)^{-8} \|u\|_{L^pL^q(Q_\rho)}^4 + C (\rho-r)^{-2}\|u\|_{L^pL^q(Q_\rho)}^2 + \frac14 \|\nabla u\|_{L^2L^2(Q_\rho)}^2.
\een

On the other hand, it follows from $\eqref{ine:estimate p h}$, $\eqref{ine:v_B u}$ and $\eqref{prop:harmonic}$, there holds
\ben\label{eq: u estimate} \nonumber
\|u \phi^\beta\|_{L^p L^\frac{6p}{3p-4}(Q_\rho)}^2 &\leq& \|v_B \phi^\beta\|_{L^p L^\frac{6p}{3p-4}(Q_\rho)}^2 + \|\nabla \pi_{h,B} \phi^\beta\|_{L^p L^\frac{6p}{3p-4}(Q_{\sigma_2})}^2 \\ \nonumber
&\leq& \|v_B \phi^\beta\|_{L^p L^\frac{6p}{3p-4}(Q_\rho)}^2 + C (\rho-\varrho)^{-\frac 3q}\|\nabla \pi_{h,B}\|_{L^p L^q(Q_{\rho})}^2 \\
&\leq& \|v_B \phi^\beta\|_{L^p L^\frac{6p}{3p-4}(Q_\rho)}^2 + C (\rho-\varrho)^{-\frac 3q}\|u\|_{L^p L^q(Q_{\rho})}^2,
\een
and
\ben\label{eq: nabla u estimate} \nonumber
\|\nabla u \phi^\beta\|^2_{L^2 L^2(Q_\rho)} &\leq& \|\nabla v_B \phi^\beta\|^2_{L^2 L^2(Q_\rho)} + \|\nabla \nabla \pi_{h,B} \phi^\beta\|^2_{L^2 L^2(Q_{\sigma_2})}\\ \nonumber
&\leq& \|\nabla (v_B \phi^\beta)\|^2_{L^2 L^2(Q_\rho)} + \|v_B \nabla \phi^\beta\|^2_{L^2 L^2(Q_\rho)} \\ \nonumber
&&+ C (\rho-\varrho)^{-\frac 3q-1}\|\nabla \pi_{h,B}\|_{L^p L^q(Q_{\rho})}^2 \\\nonumber
&\leq& \|\nabla (v_B \phi^\beta)\|^2_{L^2 L^2(Q_\rho)} + C (\rho-\varrho)^{-2} \|u\|^2_{L^p L^q(Q_\rho)} \\
&&+ C (\rho-\varrho)^{-\frac 3q-1}\|u\|_{L^p L^q(Q_{\rho})}^2,
\een
Combining $\eqref{ine:energy'}$, $\eqref{eq: u estimate}$ and $\eqref{eq: nabla u estimate}$, we arrive at
\beno
&&\|u \phi^\beta\|_{L^p L^\frac{6p}{3p-4}(Q_\rho)}^2 +\|\nabla u \phi^\beta\|^2_{L^2 L^2(Q_\rho)} \\
&\leq& C (\rho-r)^{-\frac2{\alpha-1}} \|u\|_{L^pL^q(Q_\rho)}^{\frac{2\alpha}{\alpha-1}} + C (\rho - r)^{-8} \|u\|_{L^pL^q(Q_\rho)}^4 \\
&&+ C (\rho-r)^{-2}\|u\|_{L^pL^q(Q_\rho)}^2 + \frac34 \|\nabla u\|_{L^2L^2(Q_\rho)}^2.
\eeno
Finally, using the  iterative lmma  (see \cite[Lemma V.3.1,   p.161 ]{[Giaquinta]}),  the following Caccioppoli inequality holds
\beno
\|u\|^2_{L^pL^{\frac{6p}{3p-4}}(Q_\frac34)}+\|\nabla u\|^2_{L^2L^2(Q_\frac34)}\leq C \|u\|^2_{L^pL^q(Q_{1})} + C \|u\|^4_{L^pL^q(Q_{1})} + C \|u\|_{L^pL^q(Q_1)}^{\frac{2\alpha}{\alpha-1}}.
\eeno
for $2\leq p\leq 3$. The proof is complete.

\section{Proof of Theorem \ref{onescale}}

\subsection{\bf Case I: $2 \leq p<3$} 
The proof is divided into three steps.

{\bf Step I: Decay estimates from local energy inequality.}

Choose a cut-off function as in \cite{CKN}.
Let $G(x,t) = (4\pi t)^{-\frac32} \exp{(-\frac{|x|^2}{4t})}$ be the Gaussian kernel. For $r>0$, denote
\beno
\Phi(x,t) = r^2 G(x,r^2-t), \quad (x,t) \in \mathbb{R}^3 \times (-\infty,0),
\eeno
By direct calculation, there holds for any $0 < 4r < \rho \leq \frac12$,
\ben\label{ine:psi 123} \nonumber
&&\Phi(x,t) \geq C r^{-1}, \quad (x,t) \in Q_r;\\
&&\Phi(x,t) \leq C r^{-1},\quad |\nabla \Phi(x,t)| \leq C r^{-2}, \quad (x,t) \in Q_\rho;\\ \nonumber
&&\Phi(x,t) \leq C r^2 \rho^{-3}, \quad |\nabla \Phi(x,t)| \leq C r^2 \rho^{-4}, \quad (x,t) \in Q_\rho \setminus Q_{\frac\rho2}.
\een
Let $\eta:\mathbb{R}^3 \times \mathbb{R} \rightarrow [0,1]$ be a smooth cut-off function suitable on $Q_\rho \setminus Q_{\frac\rho2}$ with $|\partial_t \eta| + |\nabla^2 \eta| \leq C \rho^{-2}$ and $|\nabla \eta| \leq C \rho^{-1}$. Substitute $\phi = \Phi \eta$ in the local energy inequality.
Obviously,
\beno
\partial_t \phi + \Delta \phi = (\partial_t \Phi + \Delta \Phi) \eta + \Phi \partial_t \eta + 2 \nabla \Phi \cdot \nabla \eta + \Phi \Delta \eta.
\eeno
Noting that $\partial_t \Phi + \Delta \Phi = 0$, we arrive at
\ben\label{ine:phi 2}
|\partial_t \phi + \Delta \phi| \leq C r^2 \rho^{-5}.
\een
Similarly,
\ben\label{ine:phi 1}
|\nabla \phi| = |\nabla \phi \eta + \phi \nabla \eta| \leq C r^{-2} + C r^2 \rho^{-4} \leq C r^{-2}.
\een

Take a fixed ball $B=B_{\frac34}$ for $v_B$. Write $v=v_B=u+\nabla \pi_h$, $\nabla \pi_h = -E_{B_\frac34}(u)$, $\nabla \pi_1 = - E_{B_\frac34}(u \cdot \nabla u)$ and $\nabla \pi_2 = E_{B_\frac34}(\Delta u)$. Then it follows from the local energy inequality (\ref{eq:local energy inq}):
\ben\label{loce new} \nonumber
&&\int_{B_\frac34}|v(x,s)|^2\phi(x,s)dx+2\int_{Q_\frac34}|\nabla v(x,\tau)|^2\phi(x,\tau) dxd\tau\\ \nonumber
&\leq& \int_{Q_\frac34}|v(x,\tau)|^2(\partial_t\phi+\Delta \phi)dxd\tau+\int_{Q_\frac34}|v_B|^2(v-\nabla \pi_h)\cdot\nabla \phi dxd\tau\\ \nonumber
&&+2\int_{Q_\frac34} v_B\cdot\nabla \nabla \pi_h v \phi dxd\tau-2\int_{Q_\frac34} \nabla \pi_h\cdot\nabla \nabla \pi_h v \phi dxd\tau\\
&&+2\int_{Q_\frac34} (\pi_1+\pi_2)v_B\cdot \nabla\phi dxd\tau,
\een
Using $\eqref{ine:psi 123}_1$, $\eqref{ine:phi 2}$ and $\eqref{ine:phi 1}$, we have
\beno
&&\frac 1 r\|v\|^2_{L^\infty L^2(Q_r)}+\frac 1 r\|\nabla v\|^2_{L^2 L^2(Q_r)} \\
&\leq& C\left(\frac{r}{\rho}\right)^2\frac{1}{\rho^3} \int_{Q_\rho}|v|^2dxdt+ C \frac{1}{r^2} \int_{Q_\rho}|v|^3dxdt+C\frac{1}{r^2}\int_{Q_\rho} |v|^2|\nabla \pi_h|dxdt\\
&&+ C \frac{1}{r}\int_{Q_\rho} |v|^2|\nabla \nabla \pi_h|dxdt+C\frac{1}{r}\int_{Q_\rho} |\nabla \pi_h||\nabla \nabla \pi_h||v|dxdt\\
&&+ C \frac{1}{r^2}\int_{Q_\rho} |v \pi_2|dxdt+C\frac{1}{r^2}\int_{Q_\rho} |v  (\pi_1-(\pi_1)_{B_\rho})|dxdt\\
&:=& M_1 + M_2 + \cdots + M_7.
\eeno
Denote
\beno
I(r) = r^{-1} \|v\|_{L^\infty L^2(Q_r)}^2 + r^{-1} \|\nabla v\|_{L^2L^2(Q_r)}^2 + r^{-\frac94} \|\pi_1-(\pi_1)_{B_r}\|_{L^1 L^2(Q_r)}^\frac32,
\eeno
with $(\pi_1)_{B_r} = |B_r|^{-1} \int_{B_r} \pi_1$. Then we estimate $M_1-M_7$ term by term.

{\bf Step II: Growth estimate of $I(r)$.}


{\bf \underline{Estimate of $M_1$:}}
\ben\label{eq:M1}
M_1 \leq C\left(\frac{r}{\rho}\right)^2 \rho^{-1} \sup_{t \in (-\rho^2,0)}\int_{B_\rho} |v|^2 dx \leq C \left(\frac{r}{\rho}\right)^2 I(\rho).
\een

{\bf \underline{Estimate of $M_2$:}}
Using H\"{o}lder's inequality and embedding theorem, we have
\ben\label{eq:M2}
M_2 \leq C \left(\frac{\rho}{r}\right)^2 \rho^{-\frac32} \|v\|_{L^\frac{10}3L^\frac{10}3(Q_\rho)}^3 \leq C \left(\frac{\rho}{r}\right)^2 I(\rho)^\frac32.
\een

{\bf \underline{Estimate of $M_3$:}}
Using $\eqref{prop:harmonic}$  in $B_{\frac34}$ by noting that $\rho<\frac12$, it follows from H\"{o}lder's inequality that
\ben\label{eq:M3} \nonumber
M_3 &\leq& C r^{-2} \|v\|_{L^\infty L^2(Q_\rho)} \|v\|_{L^2 L^6(Q_\rho)} \|\nabla \pi_h\|_{L^2 L^3(Q_\rho)} \\ \nonumber
&\leq& C \left(\frac{\rho}{r}\right)^2 \left(\rho^{-1}\|v\|_{L^\infty L^2(Q_\rho)}^2 + \rho^{-1} \|\nabla v\|_{L^2L^2(Q_\rho)}^2 \right) \rho^{-1} \|\nabla \pi_h\|_{L^2 L^3(Q_\rho)} \\ \nonumber
& \leq& C \left(\frac{\rho}{r}\right)^2 I(\rho) \|\nabla \pi_h\|_{L^2 L^q(Q_\frac34)} \\
& \leq& C \|u\|_{L^p L^q(Q_1)} \left(\frac{\rho}{r}\right)^2 I(\rho).
\een

{\bf \underline{Estimate of $M_4$:}}
Similarity,  $M_4$ can be controlled by
\ben\label{eq:M4}
M_4 \leq C \|u\|_{L^p L^q(Q_1)} \left(\frac{\rho}{r}\right)^2 I(\rho).
\een

{\bf \underline{Estimate of $M_5$:}}
For the term $M_5$, the parts of $\nabla \pi_h$ could be controlled by $u$. Then by H\"{o}lder's inequality, using  $\eqref{prop:harmonic}$ again, we have
\ben\label{eq:M5} \nonumber
M_5 &\leq& C r^{-1} \|v\|_{L^\infty L^2(Q_\rho)} \|\nabla \pi_h\|_{L^2 L^6(Q_\rho)} \|\nabla \nabla \pi_h\|_{L^2L^3(Q_\rho)} \\ \nonumber
&\leq& C \frac \rho r I(\rho)^\frac12 \|\nabla \pi_h\|_{L^2L^q(Q_\frac34)} \|\nabla \nabla \pi_h\|_{L^2L^3(Q_\rho)} \\ \nonumber
&\leq& C \frac \rho r I(\rho)^\frac12 \|\nabla \pi_h\|_{L^2L^q(Q_\frac34)} \rho \|\nabla \pi_h\|_{L^2L^q(Q_\frac34)} \\
&\leq& \delta I(\rho) + C(\delta) \left(\frac{\rho}{r}\right)^2 \|u\|_{L^pL^q(Q_1)}^4,
\een
where $\delta>0$, to be decided, and we used $p\geq 2.$

{\bf \underline{Estimate of $M_6$:}}
Noting that $\nabla \pi_2 = -E_{B_\frac34}(\Delta u)$ is harmonic and $\|\pi_2\|_{L^2(B_{\frac34})} \leq C \|\nabla u\|_{L^2(B_\frac34)}$ by $\eqref{ine:estimate p 2}$, by $\eqref{prop:harmonic}$ there holds
\ben\label{eq:M6} \nonumber
M_6 &\leq& C r^{-2} \|v\|_{L^\infty L^2(Q_\rho)} \|\pi_2\|_{L^1 L^2(Q_\rho)} \leq C \left(\frac{\rho}{r}\right)^2 I(\rho)^\frac12 \|\pi_2\|_{L^1L^2(Q_\frac34)} \\
&\leq& C \left(\frac{\rho}{r}\right)^2 I(\rho)^\frac12 \|\nabla u\|_{L^1 L^2(Q_\frac34)} \leq \delta I(\rho) + C(\delta) \left(\frac{\rho}{r}\right)^4 \|\nabla u\|_{L^2L^2(Q_\frac34)}^2.
\een

{\bf \underline{Estimate of $M_7$:}}
The H\"{o}lder's inequality, Young's inequality and $\eqref{ine:estimate p 1}$ imply that
\ben\label{eq:M7}\nonumber
M_7 &\leq& C r^{-2} \|v\|_{L^\infty L^2(Q_\rho)} \|\pi_1 - (\pi_1)_{B_\rho}\|_{L^1 L^2(Q_\rho)} \\ \nonumber
&\leq& C(\delta) \left(\frac{\rho}{r}\right)^6 I(\rho)^\frac32 + \delta \left(\rho^{-\frac32} \|\pi_1-(\pi_1)_{B_\rho}\|_{L^1 L^2(Q_\rho)}\right)^\frac32 \\
&\leq& C(\delta) \left(\frac{\rho}{r}\right)^6 I(\rho)^\frac32 + \delta I(\rho).
\een

{\bf \underline{Estimate of $r^{-\frac32} \|\pi_1-(\pi_1)_{B_r}\|_{L^1 L^2(Q_r)}$:}}
Noting that the function $\pi_1$ satisfies
\beno
-\Delta v_1 + \nabla \pi_1 = -u \cdot \nabla u, \quad \nabla \cdot v_1 = 0 \quad {\rm in} \quad B_\frac34,
\eeno
we have $-\Delta \pi_1 = \partial_i \partial_j (u_i u_j)$ in $B_\frac34$.
Let $\zeta$ be a cutoff function which equals $1$ in $Q_{\frac\rho2}$ and vanishes outside of $Q_\rho$ with $0<4r<\rho<\frac12$. Set $\pi_1-(\pi_1)_{B_r} = p_1-(p_1)_{B_r} + p_2-(p_2)_{B_r}$ with
\beno
p_1 = \frac1{4\pi} \int_{\mathbb{R}^3} \frac1{|x-y|} \partial_i \partial_j (u_iu_j\zeta)(y) dy,
\eeno
and $p_2-(p_2)_{B_r}$ is harmonic function in $Q_{\frac\rho2}$. For any $p'>1$, according to the Calder\'{o}n-Zygmund inequality,  we deduce that
\ben\label{ine:p1}
\int_{B_\frac\rho2} |p_1|^{p'} dx \leq C \int_{B_\rho} |u|^{2p'} dx.
\een
And for the harmonic part, there holds
\ben\label{ine:p2} \nonumber
&&\int_{B_r} |p_2-(p_2)_{B_r}|^{p'} dx \leq C \left(\frac r\rho\right)^{3+p'}\int_{B_{\frac\rho2}} |p_2-(p_2)_{B_\frac\rho2}|^{p'}dx \\ \nonumber
&\leq& C \left(\frac r\rho\right)^{3+p'}\int_{B_{\frac\rho2}} |p_1-(p_1)_{B_\frac\rho2}|^{p'}dx + C \left(\frac r\rho\right)^{3+p'}\int_{B_{\frac\rho2}} |\pi_1-(\pi_1)_{B_\frac\rho2}|^{p'}dx \\
&\leq& C \left(\frac r\rho\right)^{3+p'} \int_{B_\rho} |u|^{2p'} dx + C \left(\frac r\rho\right)^{3+p'}\int_{B_{\frac\rho2}} |\pi_1-(\pi_1)_{B_\frac\rho2}|^{p'}dx.
\een
Specially, for $p'=2$, $\eqref{ine:p1}$ and $\eqref{ine:p2}$ imply that
\beno
\|\pi_1-(\pi_1)_{B_r}\|_{L^2(B_r)} \leq C \|u\|_{L^4(B_\rho)}^2 + C \left(\frac r\rho\right)^\frac52 \|\pi_1-(\pi_1)_{B_\rho}\|_{L^2(B_\rho)}.
\eeno
Integration the above inequality with respect to $t$ from $-r^2$ to $0$, there holds
\beno
\|\pi_1-(\pi_1)_{B_r}\|_{L^1 L^2(Q_r)}\leq C \|u\|_{L^2 L^4(Q_\rho)}^2 + C \left(\frac r\rho\right)^\frac52 \|\pi_1-(\pi_1)_{B_\rho}\|_{L^1 L^2(Q_\rho)}.
\eeno
Using H\"{o}lder's inequality, $\eqref{ine:v_B u}$ and $\eqref{prop:harmonic}$, we have
\beno
r^{-\frac32} \|\pi_1-(\pi_1)_{B_r}\|_{L^1 L^2(Q_r)}
&\leq& C \left(\frac\rho r\right)^\frac32 \rho^{-\frac32} \|\nabla \pi_h\|_{L^2 L^4(Q_\rho)}^2 \\
&&\quad + C \left(\frac\rho r\right)^\frac32 \rho^{-\frac32} \|v\|_{L^2 L^4(Q_\rho)}^2 + C\frac r\rho I(\rho)^\frac23 \\
&\leq& C\left(\frac\rho r\right)^\frac32 \|u\|_{L^pL^q(Q_1)}^2 + C \left(\frac\rho r\right)^\frac32 I(\rho) + C\frac r\rho I(\rho)^\frac23,
\eeno
which means that
\ben\label{ine:I 3}
r^{-\frac94} \|\pi_1-(\pi_1)_{B_r}\|_{L^1 L^2(Q_r)}^\frac32 \leq C \left(\frac \rho r\right)^{\frac94} I(\rho)^\frac32 + C\left(\frac r\rho\right)^\frac32 I(\rho) + C \left(\frac \rho r\right)^\frac94 \|u\|_{L^pL^q(Q_1)}^3.
\een

Collecting $\eqref{eq:M1} - \eqref{eq:M7}$ and $\eqref{ine:I 3}$, there holds
\ben\label{ine:I r I rho} \nonumber
I(r) &\leq& C \left(\frac r\rho\right)^\frac32 I(\rho) + C \left(\frac \rho r\right)^6 I(\rho)^\frac32 + C \|u\|_{L^pL^q(Q_1)} \left(\frac \rho r\right)^2 I(\rho) \\ \nonumber
&&+ 3 \delta I(\rho) + C \left(\frac \rho r\right)^2 \|u\|_{L^pL^q(Q_1)}^4 + C \left(\frac \rho r\right)^4 \|\nabla u\|_{L^2(Q_\frac34)}^2 \\
&&+ C  \left(\frac \rho r \right)^\frac94 \|u\|_{L^pL^q(Q_1)}^3.
\een

{\bf Step III: Iterative arguments.}
Letting $r = \theta \rho$ for any $\theta \in (0,\frac14)$ in $\eqref{ine:I r I rho}$, there holds
\beno
I(r) &\leq& C \left(\theta^{\frac 32}+ \theta^{-2} \|u\|_{L^pL^q(Q_1)} + 3\delta\right) I(\rho) + C \theta^{-6} I(\rho)^\frac32 \\
&& + C \theta^{-2} \|u\|_{L^p L^q(Q_1)}^4 + C\theta^{-4} \|\nabla u\|_{L^2 L^2(Q_\frac34)}^2 +  C \theta^{-\frac94} \|u\|_{L^p L^q(Q_1)}^3
\eeno
Using Proposition \ref{cacci}, there holds
\beno
I(r) &\leq& C \left(\theta^{\frac 32}+ \theta^{-2} \|u\|_{L^pL^q(Q_1)} + 3\delta\right) I(\rho) + C \theta^{-6} I(\rho)^\frac32 + C \theta^{-2} \|u\|_{L^p L^q(Q_1)}^4\\
&& \quad + C\theta^{-4} \left(\|u\|_{L^p L^q(Q_1)}^2 + \|u\|_{L^p L^q(Q_1)}^4 + \|u\|_{L^p L^q(Q_1)}^\frac{2\alpha}{\alpha-1}\right) +  C \theta^{-\frac94} \|u\|_{L^p L^q(Q_1)}^3.
\eeno
Choose $\delta = \delta_0$ and $\theta = \theta_0 \in (0,\frac14)$ small enough but fixed from now on, and take $\varepsilon_0$ is small such that
\beno
C\theta^\frac32 + C \theta^{-2} \|u\|_{L^pL^q(Q_1)} + 3\delta \leq C\theta_0^\frac32 + C \theta_0^{-2} \varepsilon_0 + 3\delta_0 \leq \frac14,
\eeno
and
\beno
C \theta_0^{-2} \varepsilon_0^4 + C \theta_0^{-\frac94} \varepsilon_0^3 + C \theta_0^{-4} (\varepsilon_0^2 + \varepsilon_0^4 + \varepsilon_0^\frac{2\alpha}{\alpha-1}) \leq C \varepsilon_0^2.
\eeno
Then we arrive at
\ben\label{ine:iteration}
I(\theta_0 \rho) = I(r) \leq \frac14 I(\rho) + C I(\rho)^\frac32 + C \varepsilon_0^2.
\een

Similar estimates as  (\ref{ine:energy''}) by scaling or choosing a different domain (for example,  for $B=B_{\frac34}$), we have
\ben\label{ine:energy'''}\nonumber
&&\|v\phi^\beta\|^2_{L^\infty L^2(Q_\frac34)}+\|\nabla (v\phi^\beta)\|^2_{L^2 L^2(Q_\frac34)} \leq C  \|u\|_{L^pL^q(Q_\frac34)}^{\frac{2\alpha}{\alpha-1}} + C \|u\|_{L^pL^q(Q_\frac34)}^4 \\
&&
+ C \|u\|_{L^pL^q(Q_\frac34)}^2 + \frac14 \|\nabla u\|_{L^2L^2(Q_\frac34)}^2.
\een
and at this time $\phi=1$ in $Q_{1/2}$.
Then
\beno
\|v\|^2_{L^\infty L^2(Q_\frac14)}+\|\nabla v\|^2_{L^2 L^2(Q_\frac14)} &\leq& C \|u\|_{L^pL^q(Q_1)}^{\frac{2\alpha}{\alpha-1}} + C \|u\|_{L^pL^q(Q_1)}^4 \\ \nonumber
&& + C\|u\|_{L^pL^q(Q_1)}^2 + \frac14 \|\nabla u\|_{L^2L^2(Q_{\frac34})}^2.
\eeno
Using Propositon \ref{cacci} agian, there holds
\ben\label{ine:v}
\|v\|^2_{L^\infty L^2(Q_\frac14)}+\|\nabla v\|^2_{L^2 L^2(Q_\frac14)} \leq C \left(\varepsilon_0^2 + \varepsilon_0^{\frac{2\alpha}{\alpha-1}} + \varepsilon_0^4\right).
\een
For the pressure part of $I(\frac14)$,
if $q \geq 4$,  noting that $\|\pi_1-\bar{\pi}_1\|_{L^1L^2(Q_\frac14)} \leq C \|u\|_{L^2L^4(Q_\frac34)}^2$ due to (\ref{ine:estimate p 1}), there holds
\ben\label{ine:q > 4}
\|\pi_1-\bar{\pi}_1\|_{L^1L^2(Q_\frac14)} \leq C \|u\|_{L^2L^4(Q_\frac34)}^2 \leq C \|u\|_{L^pL^q(Q_\frac34)}^2.
\een
If $\frac 94 < q < 4$, by embedded inequality, we have
\ben\label{ine:q < 4}
\|\pi_1-\bar{\pi}_1\|_{L^1L^2(Q_\frac14)} \leq C \|u\|_{L^2L^4(Q_\frac34)}^2 \leq C \left(\|\nabla u\|_{L^2(Q_\frac34)}^2 + \|u\|_{L^pL^q(Q_\frac34)}^2\right).
\een
Combining $\eqref{ine:v}$, $\eqref{ine:q > 4}$ and $\eqref{ine:q < 4}$, using Proposition \ref{cacci}, there holds
\beno
I(\frac14) \leq C \left(\varepsilon_0^2 + \varepsilon_0^{\frac{2\alpha}{\alpha-1}} + \varepsilon_0^4\right) \leq \varepsilon_0^\frac32:= \varepsilon_1.
\eeno

Without loss of generality, we set $I(\rho_0)\leq \varepsilon_1$ for some $\rho_0 > 0$. Assume that for any $s \in \mathbb{N}_+$, $s < n$,
\beno
I(\theta_0^{s-1} \rho_0) \leq \varepsilon_1,
\eeno
then for $s = n$, by $\eqref{ine:iteration}$, there holds
\beno
I(\theta_0^n \rho_0) &\leq& \frac14 I(\theta_0^{n-1} \rho_0) + C I(\theta_0^{n-1} \rho_0)^\frac32 + C \varepsilon_0^2 \\
&\leq& \left(\frac14 + C \varepsilon_1^\frac12 + C \varepsilon_1^\frac13\right) \varepsilon_1.
\eeno
Choosing $\varepsilon_1$, which is dependent on $\varepsilon_0$, is small enough but fixed, such that $\frac14 + C \varepsilon_1^\frac12 + C \varepsilon_1^\frac13 \leq 1$, we arrive $I(\theta_0^n \rho_0) \leq \varepsilon_1$.
By mathematical induction, for any $n \in \mathbb{N}$,
\beno
I(\theta_0^n \rho_0) \leq \varepsilon_1.
\eeno
Then for any $r \in (0,\frac14)$, there exist constant $n_0$ such that $\theta_0^{n_0}\rho_0 < r \leq \theta_0^{n_0-1} \rho_0$. Then
\beno
&&r^{-1} \|v\|_{L^\infty L^2(Q_r)}^2 + r^{-1} \|\nabla v\|_{L^2L^2(Q_r)}^2 \\
&\leq& \theta_0^{-n_0} \rho_0^{-1} \|v\|_{L^\infty L^2(Q_{\theta_0^{n_0-1}\rho_0})}^2 + \theta_0^{-n_0} \rho_0^{-1} \|\nabla v\|_{L^\infty L^2(Q_{\theta_0^{n_0-1}\rho_0})}^2 \\
&\leq& C \theta_0 I(\theta_0^{n_0-1}\rho_0) \leq C \varepsilon_1.
\eeno
By translation invariance of Navier-Stokes equations, we obtain
\ben\label{ine:energy v}
\sup_{z_0\in Q_{1/4}}\sup_{r\in(0,1/4)}\{r^{-1} \|v\|_{L^\infty L^2(Q_r)}^2 + r^{-1} \|\nabla v\|_{L^2L^2(Q_r)}^2\}\leq C \varepsilon_1.
\een

Next we prove the regularity. By triangle inequality, there holds
\beno
r^{-2} \|u\|_{L^3(Q_r)}^3 \leq r^{-2} \|v\|_{L^3(Q_r)}^3 + r^{-2} \|\nabla \pi_h\|_{L^3(Q_r)}^3.
\eeno
It is sufficient to estimate the term $r^{-2} \|\nabla \pi_h\|_{L^3(Q_r)}^3$, since $\eqref{ine:energy v}$ implies the smallness of the term $r^{-2} \|v\|_{L^3(Q_r)}^3$, and
\beno
r^{-2} \|\nabla \pi_h\|_{L^3(Q_r)}^3 &=& r^{-2} \int_{I_r} \|\nabla \pi_h\|_{L^3(B_r)}^3 \leq r^{-2} \int_{I_r} r^3 \|\nabla \pi_h\|_{L^2(B_\frac34)}^3\\
&\leq& C r^{-2} \int_{I_\frac34} r^3 \|u\|_{L^2(B_\frac34)}^3 \leq C r \|u\|_{L^\infty L^2(Q_\frac34)}^3.
\eeno
Note that $u$ is local suitable weak solution, which means that $\|u\|_{L^\infty L^2(Q_\frac34)}^3 < + \infty$. Then there exists $r_0>0$ such that $C r_0 \|u\|_{L^\infty L^2(Q_\frac34)}^3 \leq C \varepsilon_1^\frac32$, and we have
\beno
r^{-2} \|u\|_{L^3(Q_{r})}^3 \leq C \varepsilon_1^\frac32,\quad \forall~0<r\leq r_0,
\eeno
which means that for all $z_0 \in Q_\frac{r_0}2$ by Wolf's result (\ref{eq: regularity on u}) or Wang-Wu-Zhou's result (\ref{eq:wwz}).

\subsection{\bf Case II: $p \geq 3$.}

At this time, there holds $\frac32<q\leq 9$, since $1\leq \frac2p+\frac3q<2$.

{\bf \underline{Case of  $3 \leq  q \leq  9$.}} It follows that
\beno
\|u\|_{L^3 L^3(Q_\frac12)} \leq C \|u\|_{L^p L^q(Q_\frac12)}
\eeno
which implies the regularity due to (\ref{eq: regularity on u}) or (\ref{eq:wwz}).

{\bf \underline{Case of  $\frac 32 < q< 3$.}} It follows from Proposition \ref{cacci} that
\beno
 \|u\|_{L^p L^{\frac{6p}{3p-4}}(Q_\frac12)}+\|\nabla u\|_{L^2 L^2(Q_\frac12)}\leq C\left(\|u\|_{L^pL^q(Q_1)} + \|u\|_{L^pL^q(Q_1)}^2 + \|u\|_{L^pL^q(Q_1)}^{\frac{\alpha}{\alpha-1}}\right)
\eeno
where $\alpha = \frac 2{\frac2p+\frac3q}$.
Thus
\beno
\|u\|_{L^3 L^3(Q_\frac12)} \leq C\|u\|_{L^3 L^\frac{18}{5}(Q_\frac12)} \leq C \left(\|u\|_{L^p L^q(Q_\frac12)}+\|\nabla u\|_{L^2 L^2(Q_\frac12)}\right).
\eeno
Apply Wolf's result again. The proof is complete.

\noindent {\bf Acknowledgments.}
 W. Wang was supported by NSFC under grant 12071054,  National Support Program for Young Top-Notch Talents and by Dalian High-level Talent Innovation Project (Grant 2020RD09).
 D. Zhou was supported by NSFC under grant 12071113.

\bibliographystyle{unsrt}

\end{document}